\documentclass[leqno,12pt]{amsart}
\usepackage{amssymb,amsmath,amscd}

\newtheorem{thm}{Theorem}[section]
\newtheorem{monodromy}[thm]{Monodromy Theorem}
\newtheorem{corollary}[thm]{Corollary}
\newtheorem{prop}[thm]{Proposition}
\newtheorem{lemma}[thm]{Lemma}
\newtheorem{fact}[thm]{Fact}

\theoremstyle{definition}
\newtheorem{defn}[thm]{Definition}
\newtheorem{example}[thm]{Example}

\theoremstyle{remark}
\newtheorem{remark}[thm]{Remark}

\newcommand{\bt}{\begin{thm}}
\newcommand{\et}{\end{thm}}
\newcommand{\bp}{\begin{prop}}
\newcommand{\ep}{\end{prop}}
\newcommand{\bd}{\begin{defn}}
\newcommand{\ed}{\end{defn}}
\newcommand{\bl}{\begin{lemma}}
\newcommand{\el}{\end{lemma}}
\newcommand{\bfa}{\begin{fact}}
\newcommand{\efa}{\end{fact}}
\newcommand{\bc}{\begin{corollary}}
\newcommand{\ec}{\end{corollary}}
\newcommand{\bex}{\begin{example}}
\newcommand{\eex}{\end{example}}
\newcommand{\br}{\begin{remark}}
\newcommand{\er}{\end{remark}}
\newcommand{\ben}{\begin{enumerate}}
\newcommand{\een}{\end{enumerate}}

\newcommand{\sotto}[2]{#1_{#2}}

\newcommand{\ra}{\rightarrow}

\newcommand{\ideal}[1]{\sotto {{\mathcal I}}{#1}}

\newcommand{\pso}{\mathbb{P}^3}
\newcommand{\pone}{\mathbb{P}^1}
\newcommand{\PP}{\mathbb{P}}

\newcommand{\coo}{{\mathcal O}}

\newcommand{\C}{\widetilde{C}}

\begin{document}

\title{Monodromy of Projective Curves}

\author{Gian Pietro Pirola and Enrico Schlesinger}

\address{Dipartimento di Matematica "F. Casorati", Universit\'{a} di Pavia,
via Ferrata 1, 27100 Pavia,  Italia}

\email{pirola@dimat.unipv.it}

\address{Dipartimento di Matematica, Politecnico di Milano, Piazza Leonardo da
Vinci 32, 20133 Milano, Italia}

\email{enrsch@mate.polimi.it}

\thanks{The first author was partially supported by:
1) MIUR PRIN 2003:
  {\em Spazi di moduli e teoria di Lie};
  2) Gnsaga; 3) Far 2002 (PV):
  {\em Variet\`{a} algebriche, calcolo algebrico, grafi orientati e
topologici.}
The second author was partially supported by
MIUR PRIN 2002 {\em Geometria e classificazione delle variet\`a proiettive complesse}.
}

\begin{abstract}
The uniform position principle states that, given an irreducible nondegenerate
curve $C \subset \PP^r (\mathbb{C})$,
a general $(r\!-\!2)$-plane $L \subset \PP^r$ is {\em uniform}, that is,
projection from $L$ induces a rational map $C \dashrightarrow \PP^{1}$
whose monodromy group is the full symmetric group.
In this paper we first show the locus of non-uniform $(r-2)$-planes has
codimension at least two in the Grassmannian. This result is sharp because,
if there is a point $x \in \PP^r$ such that projection from $x$ induces
a map $C \dashrightarrow \PP^{r-1}$ that is not birational onto its image,
then the Schubert cycle $\sigma(x)$ of $(r\!-\!2)$-planes through $x$
is contained in the locus of non-uniform subspaces. For a smooth curve $C$
in $\PP^3$, we show any irreducible surface of non-uniform lines is a Schubert cycle
$\sigma(x)$ as above, unless $C$ is a rational curve of degree three, four or six.
\end{abstract}

\maketitle

\section{Introduction}
The monodromy group of a branched covering of smooth complex curves
has been object of research since the early days of
algebraic geometry. Zariski~\cite{Zariski1} showed that, for a general smooth
complex projective
curve $X$ of genus $g > 6$, the monodromy group of any covering
$X \ra \PP^1$ is not solvable.
This result has been greatly generalized recently and it is now known that,
for a general curve $X$ of genus $g>3$, the monodromy group of a degree
$d$ indecomposable covering $X \ra \PP^1$ is either $S_d$,
in which case we say the covering is {\it uniform},
or $A_d$~\cite{gm}. As far as existence,
it is well known that a general (resp. every) curve $X$ of genus $g$ admits a degree $d$ uniform
covering $X \ra \PP^1$ if $d \geq \frac{g+2}{2}$ (resp. if $d \geq g+2$), while only
recently it has been established that
a general (resp. every) curve $X$ of genus $g$ admits a degree $d$
covering $X \ra \PP^1$ with monodromy group $A_d$
if $d \geq  2g+1$~\cite{mv} (resp. if $d \geq 12g+4$~\cite{artebani}).

Recently some interesting results have appeared
on the monodromy groups of coverings obtained projecting a smooth plane curve from a point.
Miura and Yoshihara~\cite{miura,yoshi-miura,yoshidual}  have determined all groups
occurring as monodromy groups of projections of smooth plane curves of degree $d \leq 5$.
Cukierman~\cite{cukierman} has shown that, for a general smooth plane curve
$C \subset \PP^2$ of degree $d$ and for every point $x \in \PP^2 \smallsetminus C$,
projection from $x$ induces a {\it uniform} covering $C \ra \PP^1$.

In this paper we  begin a systematic study of coverings obtained by projections of
curves embedded in projective space.
Our starting point is the classical {\it uniform position principle}.
Consider a degree $d$ irreducible algebraic curve
$C \subset \PP^r$ that is not contained in a hyperplane.
The uniform position principle, in the formulation
of Joe Harris~\cite{harris-eis}, states
that, if $U$ is the set of hyperplanes meeting $C$ transversally,
then the monodromy group of the covering
\begin{equation} \label{harrispi}
 \pi:  \{ (x,H) \in C \times U: \;\; x \in H \}
\ra U\, ,  \,\,\,\,\;\;\; \pi(x,H) = H,
\end{equation}
is the full symmetric group $S_d$.
This statement is equivalent - by a Lefschetz type theorem on
the fundamental group due to Zariski~\cite{Zariski} -
to the assertion that a {\it general}  $(r\!-\!2)$-plane  $L \subset \PP^r$
is uniform, that is, projection from $L$ exhibits
the normalization $\C$ of $C$ as a uniform covering of $\PP^1$
(cf. Proposition~\ref{uniform}).

Our first result is
\bt \label{cotwo}
In the Grassmannian $\mathbb{G}(r\!-\!2, \PP^r)$
the locus of non-uniform subspaces
has codimension at least two.
\et
Note that this theorem strengthen the uniform position principle.
For smooth curves it is an easy result (cf. Remark~\ref{smooth-two}), but for
curves with arbitrary singularities the proof is more involved
and is carried out in Section~\ref{three}, where we also consider subspaces
that meet $C$. In particular, for a curve in $\PP^2$, our result says
only finitely many points of $\PP^2$, including those lying on $C$,
may fail to be uniform.
The theorem is sharp because there are curves with
codimension two families of non-uniform subspaces.
To see this,
observe a subspace $L$ is certainly not uniform if it is
{\it decomposable}, that is, the projection $\pi_L: \C \ra \PP^1$
factors as the composition of two morphisms
of degree at least two. And this is the case if $L$ contains
a non birational point $x$, that is, a point $x$ such that projection from
$x$ induces a morphism $\pi_x: \C \ra \PP^{r-1}$ that is not birational
onto its image. Thus the Schubert cycle $\sigma(x)$
of $(r\!-\!2)$-planes through a non birational point $x$  is contained in
the non-uniform locus.

In Section~\ref{four} we give strong restrictions on families of decomposable
$(r\!-\!2)$-planes. For a non rational curve $C \subset \PP^3$, we show
any codimension two irreducible family of decomposable lines
is in fact a Schubert cycle $\sigma(x)$ for some non birational point $x \in \PP^3$.
More generally, we prove:
\bt
Suppose $C \subset \PP^r$ is a nondegenerate irreducible curve,
and $\Sigma \!\subset\! \mathbb{G} (r\!-\!2, \PP^r)$ is an
irreducible closed subscheme of dimension $r\!-\!1$ such that
\begin{enumerate}
\item
the general subspace $L \in \Sigma$ does not meet $C$ and is decomposable;
\item
every hyperplane $H$ in $\PP^r$ contains at least one subspace $L \in \Sigma$.
\end{enumerate}
Then $C$ is rational.
\et

In the case of {\it smooth} curves $C$ in $\PP^3$, we can actually
classify pairs $(C,\Sigma)$
where $\Sigma$ is a surface of non-uniform lines for $C$. The result
we prove in Sections~\ref{five} and~\ref{six} is

\bt Suppose
$C\subset\PP^3$ is a smooth irreducible nondegenerate curve,
$\Sigma \subset \mathbb{G} (1, \PP^3)$ is an irreducible surface,
and the general line in $\Sigma$ does not meet $C$ and is
not uniform.
Then one of the following three possibilities holds:
\ben
\item either there exists a non birational point
$x \in \pso$ such that $\Sigma$ is the Schubert cycle $\sigma(x)$ of
lines through $x$; or
\item
$C$ is a twisted cubic curve, and the general line in
$\Sigma$ is the intersection of two osculating planes to $C$; or
\item
$C$ is a rational curve of degree $4$ (resp. $6$), the general line $L$ in
$\Sigma$ is the intersection of two bitangent (resp. tritangent) planes to $C$,
and the projection $\pi_L$ is decomposable in the form
$
\C \ra \PP^1 \stackrel{\beta}{\ra} \PP^1
$
where $\beta$ has degree $2$.
\een
\et
In fact, every rational curve of degree $3$ or $4$ admits a surface of non-uniform
lines as in the theorem. The rational sextic curve studied in~\cite{bm}
admits a surface of non-uniform lines, but we show this is not
the case for a general rational sextic - see Theorem~\ref{end}.

We close the paper, in Section~\ref{seven}, giving examples of one dimensional families
of non-uniform lines for curves $C$ of arbitrary genus.

We work over the field $\mathbb{C}$ of complex numbers. If $V$ is a $\mathbb{C}$-vector
space, we denote by $\PP(V)$ the projective space of lines in $V$.
The sentence {\it for a general point $x$ in an algebraic variety $X$}
means for all points $x$ in some Zariski open dense subset of $X$.

We would like to thank  Irene Sabadini and Hisao Yoshihara for helpful conversation,
and Sonia Brivio who carefully read a preliminary version of the manuscript.

\section{Preliminaries} \label{prel}
We collect in this preliminary section some well known results
about branched coverings of algebraic curves. For the convenience of the reader,
we include proofs of some of these results.

We begin recalling
the definition of the {\it Galois group} or {\it monodromy group}
of a finite morphism of algebraic varieties~\cite{harrisgalois}.
Suppose $X$ and $Y$ are irreducible algebraic varieties of the
same dimension over $\mathbb{C}$, and $\pi: X \ra Y$ is a
generically finite dominant morphism of degree $d$. Then the
function field $K(X)$ is a finite algebraic extension of $K(Y)$ of
degree $d$. Let $\hat{K}$ denote a Galois closure of $K(X)/K(Y)$. We define
the Galois group $G=G_{\pi}$ of the morphism $\pi: X \ra Y$ to be the Galois
group of the extension $\hat{K}/K(Y)$:
the group $G$ consists of the automorphisms of the
field $\hat{K}$ that leave every element of $K(Y)$ fixed.  Given a general
point $q \in Y$, the group $G$ acts faithfully on the fibre
$\pi^{-1} (q)$ and thus may be regarded as a subgroup of
$\mbox{Aut} (\pi^{-1} (q)) \cong S_d$, and as such it is a transitive subgroup.
There is an equivalent, more geometric description of the Galois group. We can
choose a Zariski open dense subset $U \!\subset\! Y$ such that the
restriction of $\pi$ to $\pi^{-1} (U) \ra U$ is \'etale, thus a
$d$-sheeted covering in the classical topology. If $q$ is a point
of $U$, lifting loops at $q$ to $\pi^{-1} (U)$ we obtain the
monodromy representation $\rho: \pi_1 (U,q) \ra \mbox{Aut}
(\pi^{-1} (q))$. It is not difficult to see that the image of
$\rho$, that is, the monodromy group of the covering, is isomorphic to
the Galois group $G$, and is therefore independent of the
choice of $U$~\cite[Proposition p. 689]{harrisgalois}.
\bd
Let $\pi: X \ra Y$ be a generically finite dominant morphism of degree $d$
between complex algebraic varieties. We say that $\pi$ is {\it uniform}
if the monodromy group of $\pi$ is the full symmetric group $S_d$.
We say that $\pi$ is {\it decomposable} if there exists an open dense subset $U \subseteq Y$
over which $\pi$ factors as
$$
\pi^{-1} (U) \stackrel{\alpha}{\ra} V \stackrel{\beta}{\ra} U
$$ where $\alpha$ and $\beta$ are finite morphisms
of degree at least $2$.
\ed
\br \label{decompose}
One knows $\pi$ is decomposable if and only if there is an intermediate field
in the extension $K(X)/K(Y)$, and this is equivalent to the condition the
Galois group $G_{\pi}$ be imprimitive. Furthermore:

\begin{enumerate}
\item If $\pi$ is uniform, then it is indecomposable.
\item If  $\pi$ is indecomposable and $G_{\pi}$ contains a
transposition, then $\pi$ is uniform.
\end{enumerate}

The first statement is clear. Suppose $G_{\pi} \!\subset\! S_{d}$
contains a transposition, and let $N \subseteq G_{\pi}$
be the subgroup generated by the transpositions contained in $G_{\pi}$.
Then $N$ is a nontrivial normal  subgroup of $G$. If $N$ is not transitive,
then $G_{\pi}$ is imprimitive and $\pi$ is decomposable.
On the other hand, if $N$ is  transitive, then $\pi$ is uniform because
a transitive  subgroup of $S_{d}$ generated by transpositions must be all of
$S_{d}$ - see for example~\cite{cukierman}.
\er

Throughout the paper $C$  denotes an irreducible reduced curve of degree $d$
in $\PP^r$ - except in the Section~\ref{six} where we treat rational curves of even degree $2d$.
We assume $C$ is nondegenerate, that is, not contained in a hyperplane. We denote by $\C$
the normalization of $C$, by $f: \C \ra \PP^r$ the normalization map followed by the embedding of $C$
in $\PP^r$, and by $|V|$ the base point free $g^{r}_d$
on $\C$ corresponding to $f$.  We have $|V| = \PP(V)$ where
$V \subseteq \mbox{H}^0 (\C, f^* \coo_{\PP^r} (1))$ is the subspace
generated by the restrictions of the linear forms on $\PP^r$ to $\C$.
There is a one to one correspondence between hyperplanes $H$
in $\PP^r = \PP(V^*)$ and divisors of the linear series $|V|$, and
we let $D_H = f^* (H)$ denote the  divisor corresponding to $H$.
If $H$ is transversal to $C$, then $D_H$ may be identified with
the intersection of $C$ and $H$.

Given a linear subspace $M \subset  \PP^r$ of  codimension $s \geq
1$, projection  from $M$  is the  morphism $\pi_M:  \PP^r\!\smallsetminus\! M \ra \PP^{s-1}_M$
mapping a  point $x  \notin  M$ to
the  linear subspace generated by $x$ and $M$, where $\PP^{s-1}_M
\cong \PP^{s-1}$ denotes the  Schubert cycle of codimension $s\!-\!1$
subspaces of $\PP^r$ containing $M$.  Projection from  $M$ induces
a morphism $\pi_M:\C \ra \PP^{s-1}$.

For a codimension two subspace $L$, the morphism $\pi_L: \C \ra
\PP^1_{L}$ is finite (because $C$ is nondegenerate) of degree $d_L
= d - \deg (D_L)$, where $D_L$ is the base locus of the linear
system $\{D_H: H \supset L\}$. We let $G_L$ denote the monodromy
group of $\pi_L$. We say $L$ is  {\it uniform} (resp. {\it decomposable})
if the morphism $\pi$ is uniform (resp. decomposable).

\br Let $\mathbb{G} (r\!-\!2,\PP^r)$ the Grassmannian of
codimension two subspaces of $\PP^r$. The set of uniform subspaces
(with respect to a given curve $C$) is constructible in
$\mathbb{G} (r\!-\!2,\PP^r)$, and it contains a non-empty open
set (cf.~\ref{uniform} below).
\er

The projections $\pi_L$, as $L$ varies in the set of $(r\!-\!2)$~-planes
that do not meet $C$, are restrictions of the global projection:
$$ \pi: I = \{ (P,H) \in \C \times \PP^{r*}: \;\; f(P) \in H \}
\ra \PP^{r*},  \,\,\,\, (P,H) \mapsto H.
$$
Note that $I$ is a $\PP^{r-1}$-bundle over $\C$. In fact, from the Euler sequence
one sees $I$ is the projective bundle of {\em lines} in the fibres of
$f^* \Omega_{\PP^r} (1)$. By construction the scheme theoretic fibre of $\pi$
over $H\in \PP^{r*}$ is $f^*(H) = D_H$; thus $\pi$ is a finite flat surjective  morphism
(= a branched covering) of degree $d = \deg (C)$.

We use the well known fact~\cite[p.111]{acgh}
that the monodromy of $\pi: I \ra \PP^{r*}$ is the full symmetric
group $S_d$ to show  the general codimension two subspace
$L$ is uniform - for smooth curves this is proven for example in~\cite{cukierman}.
More precisely we have:
\bp \label{uniform}
Let $L \!\subset\! \PP^{r}$ be a codimension two subspace that does not meet $C$,
and let $B_{red}$ denote the
branch divisor of $\pi:I \ra \PP^{r*}$ taken with its reduced scheme structure.
If the line $\PP^1_L$ meets $B_{red}$ transversally, then  $L$ is uniform.
\ep
\begin{proof}
Given a codimension two subspace $L$, the codomain of the projection
$\pi_L$ is the line $\PP^1_L \!\subset\! \PP^{r*}$.  If $L$ does not meet $C$, then
$\pi^{-1} (\PP^1_L) \cong \C$, and under this isomorphism the restriction
of $\pi$ to $\pi^{-1}(\PP^1_L)$ corresponds to $\pi_L$.

The monodromy group of $\pi: I \ra \PP^{r*}$ is the full symmetric
group~\cite[Lemma p.111]{acgh}.  Since $\PP^1_L$ meets $B_{red}$
transversally, by a  Lefschetz-type theorem due to
Zariski~\cite{Zariski} (see also~\cite[p. 34]{fl}
and~\cite{voisin}), the inclusion $\PP^1_L \hookrightarrow
\PP^{r*}$ induces an epimorphism $\pi_1 (U \cap \PP^1_L ) \ra
\pi_1 (U)$, where $U = \PP^{r*} \!\smallsetminus\! B_{red}$. Hence
the monodromy of $\pi_L$ is the full symmetric group as well.
\end{proof}

\br\label{dual variety}
We will see below the branch divisor $B$ contains the dual variety $C^*$ whose points
are the hyperplanes tangent to $C$. We recall the definition of $C^*$ and the biduality theorem.
Let $W^0_C \subset \PP^r \times \PP^{r*}$ denote the set of pairs $(x,H)$ where
$x$ is a smooth point of $C$ and $H$ is a hyperplane containing the tangent line
$T_x C$, and let $W_C$ denote the closure of $W^0_C$. We say that a hyperplane $H$
is tangent to $C$ at a (not necessarily smooth) point $x$ if the pair $(x,H)$ belongs to
$W_C$. The dual variety $C^{*} \!\subset\! \PP^{r*}$ is defined as the set of hyperplanes
that are tangent to $C$, i.e., as the projection of $W_C$
in $\PP^*$. Of course, these definitions apply to any subvariety of $\PP^r$. The biduality
theorem asserts $W_C = W_{C^*}$, that is, a hyperplane $H$ is tangent to $C$ at $x$ if and
only if in the dual projective space the hyperplane $x^*$ is tangent to $C^*$ at
the point $H^*$.
\er

\br
In fact, the hypothesis $L$ does not meet
$C$ is redundant in Proposition~\ref{uniform}: one can check, for example using the biduality
theorem, that  the fact $\PP^1_L$ meets $B_{red}$ transversally implies $L$ does not meet
$C$.
\er

We now wish to describe the branch divisor $B$ of $\pi$, and the subspaces
$L$ meeting $B_{red}$ non transversally. For this we need to introduce
the {\it ramification sequence} of the linear series $|V|$ at
a point $P \in \C$ as defined in~\cite{acgh}.
Since $V$ is $r+1$-dimensional, for a given $P \in \C$ the "multiplicity at $P$ function"
$H \mapsto m_{P} (D_H)$
takes on exactly $r+1$-values, so there are nonnegative integers
$\alpha_0 (P) \leq \alpha_1 (P) \leq \cdots \leq \alpha_r (P)$ such that
$$ \{ m_{P} (D_H) : \ H^* \in \PP(V) \} =
\{\alpha_0, 1+\alpha_1, 2+\alpha_2 \ldots, r+ \alpha_r \}.$$
For each $k = 0,1, \ldots, r-1$ the set of hyperplanes $H$ with
$$ m_{P} (D_H) \geq k+1 + \alpha_{k+1}(P)$$
is a linear subspace $M^{*}_{k,P}$
of $\PP (V)$ of codimension $k+1$. Its dual $f_{k}(P)$
is a $k$-dimensional subspace of $\PP(V^*)=\PP^r$ called the osculating $k$-plane
of $C$ at $P$ - of course, $f_1(P) = T_P \,C$ is the tangent line to the branch of $C$
through $P$. The osculating $k$-plane is the unique $k$-plane $M$ such that
$$m_{P} (D_M) = k+1 +\alpha_{k+1}(P)$$
(the other $k$-planes have smaller multiplicity at $P$).

One easily verifies $H$ is tangent to $C$ at $x$ if and only if there
is $P \in \C$ such that $x=f(P)$ and $T_P \, C\!\subset H$. Thus $H \in C^*$
if and only if there is
$P \in \C$ such that $m_P (D_H) \geq 2 + \alpha_2 (P)$.

We can now  describe the branch locus $B$ of $\pi:I \ra \PP^{r*}$.
For each point $x \in C$ we let
$$\beta(x) = \sum_{\{P \in \C: f(P) =x \} } \alpha_1 (P)$$
denote the total ramification of $f$ at $x$.
\bp \label{branch-divisor}
Let $C^{*} \!\subset\! \PP^{r*}$ denote the hypersurface dual to $C$, and for each
point $x \in \PP^r$ let $\Pi_x \!\subset\! \PP^{r*}$ denote the hyperplane dual to $x$.
Let $B$ denote the branch divisor of $\pi: I \ra  \PP^{r*}$. Then in the group of Weil
divisors of $\PP^{r*}$ we have:
$$
B = C^{*} + \sum_{\{x \in C\}} \beta(x) \Pi_x.
$$
Furthermore, the multiplicity of $B$ at a point $H^* \in \PP^{r*}$ is
$$
m_{H^*} (B) = \sum_{\{P \in D_H\}} (m_{P}(D_H) -1)
$$
\ep
\begin{proof}
The ramification divisor $R$ of $\pi$ is defined as the scheme of zeros on
$I$ of the Jacobian determinant of $\pi$, and the branch divisor
$B= \pi_{*} (R)$ as the push forward of $R$ to $\PP^{r*}$. Equivalently,
$R$ is defined  by the zeroth Fitting ideal of
the relative sheaf of differential $\Omega_{I/\PP^{r*}}$~\cite{kleiman}.
This definition is functorial, as is the push forward of divisors;
it follows that, if $\PP^1_L \!\subset\! \PP^{r*}$
is a line corresponding to a codimension two subspace $L$ not meeting $C$,
the branch divisor $B_L$ of $\pi_L: \pi^{-1} (\PP^1_L) \cong \C \ra \PP^1_L$
is the scheme theoretic intersection of $B$ and $\PP^1_L$.
Since $D_H= \sum_{P \in \C} m_{P} (D_H) P$ is the fibre of $\pi_L$ over $H^{*}$
and we work over $\mathbb{C}$, the branch divisor of $\pi_L$ is
$$
B_L =
\sum_{ \{H^{*} \in \PP^1_L\} } \left(\sum_{\{P \in D_H\}} (m_{P}(D_H) -1)\right)H^*.$$
Now the multiplicity of $B$ at $H^*$ is precisely the multiplicity
of intersection at $H^*$ of $B$ with a general line through $H^*$,
hence we conclude
$$
m_{H^*} (B) = \sum_{\{P \in D_H\}} (m_{P}(D_H) -1).
$$
If $H$ is a tangent hyperplane to $C$,
then there is a point $P \in \C$
appearing with multiplicity at least two in $D_H$, thus $H$ belongs to $B$ and
$C^*$ is a component of $B$.
Furthermore, since  we are in characteristic zero, the general tangent
plane $H$ to $C$ is simply tangent, that is, $m_{H^*} (B^*)=1$.
Thus $B$ is generically smooth along $C^*$, which means that $C^*$ is a
reduced component of $B$.

On the other hand, suppose $H$ is in the branch locus but not in $C^*$.
Then in the fibre $D_H$ there is a point $P$ such that
$m_{P} (D_H)  \geq 2$ because $H$ is in the branch locus and
$m_P (D_H) = 1 + \alpha_{1}(P)$ because $H$ is not tangent to $C$.
Therefore $1 + \alpha_1 (P) \geq 2$ so that $m_P (D_K) \geq 2$ for every
hyperplane $K$ through $x=f(P)$. Hence the entire hyperplane $\Pi_x \!\subset\! \PP^{r*}$
is contained in the branch locus. To compute the multiplicity of $\Pi_x$
in $B$, we pick a hyperplane $H$ through $x$ that is not tangent to $C$ and does
not pass through any other singularity of $f$, i.e. $\Pi_{x}$
is the only component of $B$ through $H^*$. Then the desired multiplicity
is equal to:
$$
m_{H^*} (B) =  \sum_{\{P \in D_H\}} (m_{P}(D_H) -1 )=
\sum_{ \{P \in \C: f(P) =x \}} \alpha_1 (P) = \beta(x).
$$

Note: since $x \in C$, by biduality $\Pi_x$ is tangent
to $C^*$.
\end{proof}

We can now list the codimension two subspaces $L$ for which
$\PP^{1}_L$ does not meet
$B_{red}$ transversally:
\bl\label{list}
Fix a codimension two subspace $L \!\subset\! \PP^r$ not meeting $C$.
Let $\Pi_{1}, \ldots, \Pi_{n}$
be the hyperplanes of $\PP^{r*}$ contained in $B$. Then $\PP^1_L$ meets $B_{red}$
transversally unless
\begin{enumerate}
\item $L$ is a contained in a hyperplane $H$ with $H^* \in C^*_{sing}$, or
\item $L$ is a contained in a hyperplane $H$ with $H^* \in C^* \cap \Pi_{i}$, or
\item $L$ is a contained in a hyperplane $H$ with $H^* \in \Pi_{i} \cap \Pi_{j}$.
\end{enumerate}
\el
\begin{proof}
The subspace $L$ satisfies one of the
conditions (1),(2),(3) if and only if  $\PP^1_L$ meets the singular locus of $B_{red}$.
Thus we only have to check that if $\PP^1_L$ is tangent to $B_{red}$ at a smooth point
$H^*$, then $L$ meets $C$. This follows from biduality if $H^{*} \in C^{*}$, while
if $H^{*} \in  \Pi_{i}=\Pi_{x_i}$ for some $i$, then
$\PP^1_L$ is tangent to and hence contained in $\Pi_{x_i}$, so that
$x_i \in L \cap C$.
\end{proof}

We close this preliminary section recalling the relation between
the monodromy group $G_L$ and the ramification type of $\pi_L$.
\bp
\label{scambi}
Let $\tau: \C \ra \PP^1$ be a finite morphism of degree $d$, and let
$H_1,\ldots, H_n~\in~\PP^1$ be the branch points of $\tau$.
Then to each $H_i$ one can associate a permutation $\sigma_i$ in the Galois group
$G_{\tau} \subseteq S_d$ so that
\begin{enumerate}
\item
the Galois group $G_{\tau}$
is generated by any  set of $n\!-\!1$ permutations among the $\sigma_i$;
\item
if $\sum_{j=1}^{t_i} m_{ij} P_{ij}$ denotes
the scheme theoretic fibre of $\tau$ over $H_i$,
then $\sigma_i$ is a permutation
with cycle structure  $(m_{i1}, \ldots, m_{it_i})$.
Thus, if $H_i$ is a simple branch point, the permutation $\sigma_i$ is a
transposition.
\end{enumerate}
In particular, if all but one of the branch points are simple, the morphism $\pi$ is uniform.
\ep

\begin{proof}
This is classical and explained for example in~\cite[\S III.4]{miranda}.
It follows from the fact the fundamental group  of
$\PP^1 \smallsetminus \{H_1, \ldots, H_n\}$  is generated  by small
loops   $\gamma_1,  \ldots,  \gamma_n$   around  $H_1,   \ldots,  H_n$
respectively,   with   the    relation   $\gamma_1   \gamma_2   \cdots
\gamma_n=1$.
 If all the branch points but one are simple, the Galois group is generated by
 transpositions, hence must be all of $S_d$.
\end{proof}

\bex \label{twisted}
Let $C \subset \PP^3$ be a twisted cubic curve, and for every $P \in C$ let
$H(P)$ denote the osculating plane to $C$ at $P$, so that $D_{H(P)}= 3 P$.
Given two distinct points $P,Q \in C$, let $L = H(P) \cap H(Q)$. Since $H(P)$ meets
$C$ only at $P$, the line $L$ does not meet $C$, and the projection $\pi_L$
has degree $3$. By Riemann-Hurwitz, the projection $\pi_L$
is ramified only over $H(P)$ and $H(Q)$, and by Lemma~\ref{scambi} the monodromy
group of $\pi_L$ is $A_3$. Letting $P$ and $Q$ vary we find a two dimensional
irreducible family of lines $L$ that are non-uniform and indecomposable. We will later
show (Theorem~\ref{space}) there is no other smooth curve in $\PP^3$ admitting a two
dimensional family of non-uniform indecomposable lines.

If we let $P$ and $Q$ come together, the line $L$ becomes tangent to $C$,
thus $\pi_L$ is an isomorphism and $L$ is uniform according to our definition.
\eex

\section{The locus of non-uniform subspaces has codimension at least two}\label{three}
In this section we prove the locus of non-uniform subspaces has
codimension at least two in the Grassmannian $\mathbb{G} (r\!-\!2, \PP^r)$.
One easily reduces this statement to the case $r=2$, and so the problem is to show
a plane curve has only finitely many non-uniform points.

\br \label{smooth-two}
If $C$ is {\it smooth} and $\mathcal{N}(C)$ denotes the locus of codimension two subspaces
that do not meet $C$ and are not uniform,  the fact that $\mathcal{N}(C)$ has
codimension at least two is an immediate consequence of
Propositions~\ref{branch-divisor} and~\ref{scambi}.
More generally suppose $C$ is a curve for which the branch divisor $B$ of the morphism
$\pi: I \ra \PP^{r*}$ is reduced, and consider a subspace $L$ in $\mathcal{N}(C)$.
Then by Proposition~\ref{scambi} there are at least two branch points $H_1$ and $H_2$
of $\pi_L$ that are not simple.
By Proposition~\ref{branch-divisor} the hyperplanes $H_1$ and $H_2$ are singular points
of the branch divisor $B$. Hence
$\PP^1_L$ meets the singular locus $B_{sing}$ of the branch divisor in at
least two distinct points. Therefore, if  $\mathcal{N}^*$ denote the locus of lines
$\{\PP^1_L : L \in \mathcal{N}(C)\}$, we have
$$\dim \mathcal{N}(C) =  \dim \mathcal{N}^* \leq 2 \dim  B_{sing} \leq 2r-4,$$
thus $\mathcal{N}(C)$ has codimension at least two in the Grassmannian.
\er

We now generalize the above statement to curves $C$ with arbitrary singularities,
without excluding subspaces that meet $C$. The key case is that of plane curves.
When $r\!=\!2$, a codimension two subspace is a point $x \in
\PP^2$, and $\PP^1_x \!\subset\! \PP^{2*}$ is the pencil of lines
through $x$. We first consider the case $x \notin C$ (outer points
in the terminology of~\cite{yoshidual}). Proposition~\ref{uniform}
implies that $x$ is uniform, unless $x$ belongs to one of the the
finitely many lines listed in~\ref{list}. Thus, if we can show
that the general point of an arbitrary line is uniform, we can
conclude that there are only finitely many non-uniform outer
points. We will need the following lemma:
\bl \label{Galois} Let $\tau:X \ra Y$ be a surjective \'etale morphism
of smooth algebraic varieties, and let $h: Y \ra U$ be a dominant morphism.
Assume that, for the general point $q \in U$,
the fibres $X_q=(h \circ \tau)^{-1} (q)$ and $Y_q = h^{-1} (q)$ are
irreducible. Then, for a general point $q \in U$, the monodromy group
of $\tau_q: X_q \ra Y_q$ is a normal subgroup of the monodromy group
of $\tau: X \ra Y$.
\el
\begin{proof}
By generic smoothness
(cf. \cite[Lemma 1.5a)]{nori}) there is a dense Zariski open subset
$U_0 \!\subset\! U$ such that $h:h^{-1} (U_0) \ra U_0$ is a fibre bundle
in the analytic topology. Thus we may assume that $h:Y \ra U$ is a fibre
bundle in the analytic topology. Let $g: Y_q \ra Y$ be the
inclusion of the fibre over $q \in U$. The long homotopy sequence
of the fibration $h: Y \ra U$ shows that $g_* \pi_1 (Y_q)$ is a
normal subgroup of $\pi_1 (Y)$.

Now, fixed a base point $y \in Y_q$, the monodromy representations
of $\tau$ and $\tau_q$ are related by the commutative diagram:
\begin{equation}\label{phi}
\begin{CD}
\pi_1 (Y_q,y) @>{\rho(\tau_q,y)}>> \mbox{Aut} (\tau_q^{-1} (y)) \\
@VV{g_*}V @VVV \\ \pi_1 (Y,y) @>{\rho(\tau,y)}>>
\mbox{Aut}(\tau^{-1} (y))
\end{CD}
\end{equation}
where the vertical arrow on the right is an isomorphism.
Since $g_* \pi_1 (Y_q)$ is a normal subgroup of $\pi_1 (Y)$,
the diagram above shows the monodromy group $\mbox{Im}(\rho(\tau_q,y))$
of $\tau_q$ is a normal subgroup of the monodromy group
$\mbox{Im} (\rho(\tau,y))$ of $\tau$.
\end{proof}

\bp \label{monodromytrick}
Let $M$ be a line in $\PP^2$. Then all but finitely many points
$x \in M$ are uniform.
\ep
\begin{proof}
Fix a general point $q$ in $M$.
We prove that the monodromy group $G_q$
of $\pi_q$ is the full
symmetric group $S_d$
showing  it contains a transposition and it is a normal subgroup of $S_d$.
To see $G_q$ contains a transposition is immediate: it is enough to observe
a general point $q$ of $M$ lies on a simple tangent line to $C$.

To see $G_q$ is normal, we use the previous lemma. For this, we need to
consider the map $\pi: I \ra \PP^{2*}$ of Proposition~\ref{uniform}
as a family of projections indexed by the points of $M$.
To this end we fix a second line $L$, let $z_0$ be the point
of intersection of $L$ and $M$, and write
$L' = L \smallsetminus \{z_0\}$, $M'= M
\smallsetminus \{z_0\}$. Then the map $L' \times M' \ra \PP^{2*}$
sending a pair $(t,u)$ to the line joining them is an open embedding.
Let $Y$ be the open subscheme $L' \times M'$ of $\PP^{2*}$ and define
$$X= \pi^{-1} (L' \times M'), \,\,\, \tau = \pi_{|X}: X \ra Y .$$
Since $Y$ is open in $\PP^{2*}$, the monodromy group of
$\tau$ coincides with the monodromy group of $\pi$, i.e., it is the symmetric group $S_d$.
Now we restrict $Y$ to an open subset over which $\tau$ is \'etale,
and apply Lemma~\ref{Galois} to $\tau$ with
$U = M'$ to conclude that, for general
$q \in M'$, the monodromy group of $\tau_q: X_q \ra Y_q$ is normal
in $S_d$.
Now it is immediate to verify that
$X_q$ is isomorphic to an open subscheme of $\C$,
and with this identification $\tau_q: X_q \ra Y_q$ is projection from $q$.
Thus the monodromy group $G_q$ is normal in $S_d$.
This finishes the proof.
\end{proof}

\br
The same argument shows that, for a nondegenerate curve $C \subset \PP^r$,
given any hyperplane $M$ in $\PP^r$, the general codimension two subspace
$L$ contained in $M$ is uniform.
\er

We can also show, with an argument similar to the proof of~\ref{monodromytrick},
that only finitely many points of $C$ (inner points) may fail to be uniform:

\bp \label{monodromytrick2}
A general point of $C$ is uniform.
\ep
\begin{proof}
Fix a general point $q$ in $C$ and
a general line $H$.
The fibre of $\pi_q$ over $H$ is $C \cap H$ with the point $q$ removed, therefore
the monodromy group $G_q$ is a subgroup of $S_{d-1}$. We may assume
$d=\deg(C) \geq 3$. Then we can find a simple tangent to $C$ passing with multiplicity
one through $q$. It follows the monodromy of $\pi_{q}$ contains a transposition.
Thus it is enough to show the monodromy group $G_q$ is normal
in $S_{d-1}$.

To this end, fix a line $L$, and let $C'$ denote the set of
smooth points of $C$ that do not belong to $L$; similarly let $L'$ denote the
set of points of $L$ that do not lie on $C$. Then $L' \times C'$ is an open
subscheme of the incidence correspondence
$$I = \{ (P,H) \in \C \times \PP^{2*}: \;\; f(P) \in H \}$$ via
the embedding $(t,x) \mapsto (f^{-1} (x), \overline{xt})$.
The projection $\pi: L' \times C' \subset I \ra \PP^{2*}$
maps a pair $(t,x)$ to the line $\overline{xt}$, and of course we know
$\pi$ is generically finite of degree $d$ with monodromy group $S_d$.

In order to be able to apply Lemma~\ref{Galois}
we need to pull back the covering $\pi$ to $L' \times C'$. This pull back, that is
the fibred product $$T=(L' \times C') \times_{\PP^{2*}} (L' \times C'),$$
has two components: one is the diagonal, the other is isomorphic to
$$
Z = \{ (x,y) \in C' \times C': x \neq y; \, \overline{xy} \cap L \in L' \}.
$$
Let $\pi_Z: Z \ra L' \times C'$ denote the composition of the inclusion
$Z \hookrightarrow T$ with the projection $pr_1: T \ra L' \times C'$, that is,
$\pi_Z (x,y) = (\overline{xy} \cap L, y)$. Then $\pi_Z: Z \ra L' \times C'$ is
generically finite of degree $d-1$,  and, for general $q \in C$, the restriction
$\pi_q: Z_q \ra L' \times \, q$
is (modulo obvious identifications) projection from $q$. Using Lemma~\ref{Galois}
we conclude the monodromy group $G_q$ is a normal subgroup of $G_{\pi_Z}$, and we will be
done if we can show $G_{\pi_Z}=S_{d-1}$.

Choose an open dense subset $Y \subset \PP^{2*}$ such that
$$\pi_{|\pi^{-1} (Y)}: W = \pi^{-1} (Y) \ra Y$$
is \'etale, and fix a base point $w \in  W$. Then as in~\ref{Galois} we look at the
diagram:
\begin{equation}\label{phi2}
\begin{CD}
\pi_1 (W,w) @>{\rho(\pi_Z)}>> \mbox{Aut} (\pi_Z^{-1} (w)) \cong S_{d-1}  \\
@VV{\pi_*}V @VVV \\ \pi_1 (Y, \pi(w)) @>{\rho(\pi)}>>
\mbox{Aut}(\pi^{-1} (\pi(w))) \cong S_d
\end{CD}
\end{equation}
Since $W \ra Y$ is a topological covering,
the image of the vertical map $\pi_*$ is the stabilizer of
$w$ under the action of $\pi_1 (Y, \pi(w))$ on the fibre $\pi^{-1} (\pi(w))$.
The bottom row of the diagram is surjective because $G_{\pi} = S_d$.
Therefore the image of $\rho(\pi) \circ \pi_*$ is the stabilizer
of $w$ in $\mbox{Aut}(\pi^{-1} (\pi(w)))$, which is also the image of
the vertical map on the right. Therefore $\rho(\pi_Z)$ is surjective, that is,
the monodromy group of $\pi_Z$ is $S_{d-1}$, and this concludes the proof.
\end{proof}
We can now prove:

\bt \label{codimension} Let $C\!\subset\!\PP^r$ be an
irreducible non degenerate curve, $r \geq 2$. In the Grassmannian
$\mathbb{G} (r\!-\!2,\PP^r)$ the locus of non-uniform subspaces  has
codimension at least two.
\et
\begin{proof}
Suppose first $r=2$. We have to show all but finitely many points $x \in \PP^2$ are uniform
for $C$.
By Proposition~\ref{uniform}, if $x$ is non uniform,
then $x$ lies either on one of the finitely many lines listed in Lemma~\ref{list} or
on $C$. By  Propositions~\ref{monodromytrick} and~\ref{monodromytrick2}
only finitely many such points may fail to be uniform.

Now assume $r \geq 3$. Fix an irreducible nontrivial effective divisor
$\mathcal{D} \subset \mathbb{G} (r\!-\!2,\PP^r)$. We show the general
subspace $L$ in $\mathcal{D}$ is uniform reducing to the case $r=2$.
For this, let $M$ be a general codimension three subspace, so that
projection from $M$
$$
\pi_M: \C \ra \PP^2_{M}
$$
is a birational morphism of $\C$ onto $C_1 = \pi_M (\C)$. Then
projecting $C$ from a codimension
two subspace $L$ containing $M$ is the same as projecting $C_1$ from the point
$\pi_M (L)$ of $\PP^2_M$.
We have seen there are only finitely many non-uniform points in $\PP^2_M$,
thus we see there are only finitely many codimension two subspaces
$L$ that contain $M$ and are not uniform.

To finish,  we need  to observe {\it  a general codimension three
subspace $M$ of $\PP^3$ is contained in infinitely many general subspaces $L$ of
$\mathcal{D}$}. More precisely, let $\sigma(M) \cong \PP^2$ denote
the set of codimension two subspaces $L$ that contain a given $M$, and fix
a proper subvariety $\Sigma \subset \mathcal{D}$. The
set of codimension three subspaces $M$ such that
$$\dim \, \sigma(M) \cap (\mathcal{D} \smallsetminus \Sigma) > 0$$
contains a dense open subset $U \subset \mathbb{G} (r\!-\!3,\PP^r)$.
Indeed, since every nontrivial effective divisor on the
Grassmannian is ample, the dimension of $\mathcal{D} \cap \sigma(M)$ is at least
one for every codimension three subspace $M$. On the other hand,
a dimension count on the incidence variety
$\{(M,L): M \subset L\}$ shows that, for $M$ general in $\mathbb{G} (r\!-\!3,\PP^r)$,
the intersection $\Sigma \cap \sigma(M)$ is at most zero dimensional.
\end{proof}

\br \label{rmk}
The bound on the codimension of the family of non-uniform subspaces is sharp.
Indeed, there are curves $C$ for which there
exist points $x \in \PP^r$ giving a non birational projection. Fix such a point $x$.
If $L$ is a codimension two subspace containing $x$, then $\pi_L$ factors through
$\pi_x$ hence  $L$ is decomposable and in particular it is non-uniform. Thus
every codimension two subspace $L$ through $x$ is non-uniform, and the family of
such subspaces has codimension two in the Grassmannian.
\er

\section{Families of decomposable subspaces}\label{four}
In this section we show that, if $C \subset \PP^r$ is not rational,
an irreducible family $\Sigma$ of decomposable $(r\!-\!2)$-planes
either has dimension less than $r\!-\!1$ or is special in the sense
that a general hyperplane of $\PP^r$ does not
contain any subspace in $\Sigma$. More precisely:

\bt \label{trick}
Suppose $\Sigma \!\subset\! \mathbb{G} (r\!-\!2, \PP^r)$ is an
irreducible closed subscheme of dimension $r\!-\!1$ such that
\begin{enumerate}
\item
the general subspace $L \in \Sigma$ does not meet $C$ and is decomposable;
\item
every hyperplane $H$ in $\PP^r$ contains at least one subspace $L \in \Sigma$.
\end{enumerate}
Then $C$ is rational.
\et

\begin{proof}
We give now a sketch of the proof, the rest of the section being
dedicated  to filling in the details. A classical theorem of
de Franchis~\cite{defranchis} states that there are only finitely many curves
$Y$ (up to isomorphisms) for which there exists a surjective morphism $\C \ra Y$.
Since the projections
$\pi_L$ form a continuous family, it follows we can find a curve $Y$ and an integer
$e$ such that, for every $L$ in an open subset $W$ of $\Sigma$, the
projection $\pi_L$ factors through some degree $e$ morphism $\alpha: \C \ra Y$,
where $\alpha$ varies in an irreducible family of morphisms $\mathcal{M}$.
Now by hypothesis a general hyperplane $H$ contains at least one subspace $L \in W$,
hence  $C \cap H$ contains a subset of $e$ points that make up a fibre
of some morphism $\alpha \in \mathcal{M}$. By the uniform position principle,
any set of $e$ points in $C \cap H$ is a fibre of some  $\alpha \in \mathcal{M}$.

If $g(Y) \geq 1$, then the set of morphisms $\C \ra Y$ of degree $e$ is finite
(up to translations in $Y$ if $g(Y)=1$), hence we can find a single morphism
$\alpha$ such that the projections $\pi_L$ with $L \in W$ all factor through $\alpha$.
But then not every set of $e$ points in $C \cap H$ is a fibre of $\alpha$,
a contradiction. Thus $Y \cong \PP^1$.

It follows  that any two fibres of a  morphism $\alpha: \C \ra Y\cong \PP^1$ are linearly
equivalent as divisors on $\C$. Now the uniform
position principle implies that in $C \cap H$ any two disjoint subsets
of cardinality $e$ are linearly equivalent as divisors on $\C$. From this one
deduces $2P$ is linearly equivalent to $2Q$ for two general points $P,Q \in \C$,
hence $C$ is rational.
\end{proof}

We now fill in the details of the proof.
\bl \label{defr}
Suppose $T$ is an irreducible variety and
$\pi: \C \times T \ra \PP^1 \times T$ is a $T$-morphism such that
for every $t \in T$ the restriction
$\pi_t: \C \times \{t\} \ra \PP^1 \times \{t\}$
is a decomposable finite morphism of degree $d$.

Then there exist an open dense
subset $W \subset T$, a curve $Y$, an integer $e$, $2 \leq e \leq d/2$,
and an irreducible family $\mathcal{M}$ of degree $e$ morphisms from $\C$ to $Y$
with the following property:
for every $t \in W$, there exists $\alpha \in \mathcal{M}$  such that
$\pi_t$ factors through $\alpha$.
\el
\begin{proof}
A well known theorem by de
Franchis~\cite{defranchis,kani,samuel} states that,
up to isomorphisms,
there are only finitely many smooth projective curves $Y$,
for which there exists a surjective morphism $\alpha: \C \ra Y$. Let
$\mbox{Mor}_d (\C,\PP^1)$ denote the quasiprojective scheme
that parametrizes the set of finite morphisms $\C \ra \PP^1$ of degree $d$
(see~\cite[I.1.10]{kollar} or \cite{debarre}); by the universal property of
$\mbox{Mor}_d (\C,\PP^1)$
the assignment $t \mapsto \pi_t$ defines a morphism
$\psi: T \rightarrow \mbox{Mor}_{d} (\C,\PP^1)$, and by assumption the image $\psi(T)$
is contained in the set of decomposable morphisms.
Now by de Franchis Theorem the set of decomposable morphisms
in $\mbox{Mor}_d (\C,\PP^1)$ is covered by finitely many sets, images of maps
$$
\mbox{Mor}_{e_j}(\C,Y_i) \times \mbox{Mor}_{d/e_j} (Y_i,\PP^1)
\ra \mbox{Mor}_{d} (\C,\PP^1),
\;\;\;\; (\alpha,\beta)\mapsto \beta \circ \alpha.
$$
where $Y_i$ varies in a finite set of curves, and
$e_j$ in the set of integers $> 1$ that divide $d$.
Since $T$ is irreducible,
there is one of these sets that contains an open dense subset
$A \subset \psi (T)$. Then we can take $W = \psi^{-1} (A)$, and
the rest is clear.
\end{proof}

We will use the following version of the
\begin{monodromy}[\cite{acgh} pp. 111-2]
Let $C\!\subset\!\PP^r$ be a reduced irreducible curve of degree $d$,
and denote by $U \!\subset\! \PP^{r*}$ the open set of hyperplanes transverse
to $C$. Then
$$
I(d) = \{(P_1, \ldots, P_d,H): H \in U; \mbox{\em  the $P_i$'s are distinct
and $D_H = P_1+ \cdots + P_d$ }  \}
$$
is an irreducible variety.
\end{monodromy}
Note that the symmetric group $S_d$ acts on the covering projection
$p: I(d) \ra U$, and $U$ is the quotient of $I(d)$ modulo the action
of $S_d$. It follows that for any subgroup $G$ of $S_d$, the quotient
space $I(d)/G$ is an irreducible variety. In particular, we have:
\bc \label{irred}
For every positive integer $m$, let $\C_m$ denote the $m$-th symmetric
product of $\C$, which we identify with the set of effective divisors of
degree $m$ on $\C$. Then the sets
\begin{enumerate}
\item
$$
I_e= \{(E,H) \in \C_e \times U : \; E \leq D_H \}
$$
\item
$$
I_{e,f}= \{(E,F,H)\in \C_e \times \C_f \times  U : \; E + F \leq D_H \}
$$
\end{enumerate}
are irreducible.
\ec
\begin{proof}
$I_e$ is the quotient of $I(d)$ by the subgroup of $S_d$ that stabilizes the
subset $\{1,\ldots ,e \}$ of $\{1, \ldots, d \}$; while $I_{e,f}$
is the quotient of $I(d)$ by the subgroup of $S_d$ that stabilizes the
subsets $\{1,\ldots ,e \}$ and $\{e+1, \ldots, e+f \}$ of $\{1, \ldots, d \}$.
\end{proof}

\begin{proof}[Proof of Theorem~\ref{trick}]
Note that by assumption $(2)$ the morphism
$$
p:
J = \{ (L,H) \in \Sigma \times \PP^{r*}: \; L \!\subset\! H \}
\ra
\PP^{r*},
\;\;\;\; (L,H)\mapsto H
$$
is dominant and generically finite. It follows that every open dense subset
$T \subset \Sigma$ has the property that, given a general hyperplane $H$, we
can find $L$ in $T$ that is contained in $H$.

Thus under the hypotheses of the theorem
we can find $T \!\subset\! \mathbb{G} (r\!-\!2, \PP^r)$
irreducible of dimension $r\!-\!1$ such that, for every $L \in T$,
we have $L \cap C = \emptyset$, the projection $\pi_L$ is decomposable,
and the general hyperplane $H$ contains some $L \in T$.
We wish to apply Lemma~\ref{defr} . For this, we
fix a line $l \subset \PP^r$ such that the open subset
$T_l = \{ L \in T : L \cap l = \emptyset \}$
of $T$ is nonempty. Setting $\PP^1 = l$, we define
$$ \pi: \C \times T_l \ra \PP^1 \times T_l$$
mapping a pair $(P,L)$ to the pair $(\pi_L(P),L)$ where by abuse of notation
$\pi_L (P)$ is the projection of $f(P)$ from $L$ to $l$.

By Lemma~\ref{defr} we can find an open dense subset $W \subset T_l$,
a curve $Y$ and an integer $e$, $2 \leq e \leq d/2$,
such that for every $L$ in $W$
the morphism $\pi_L$ factors as $$
\C \stackrel{\alpha}{\ra} Y
\ra \PP^1
$$
where $\alpha$ belongs to some irreducible component $\mathcal{M}$
of $\mbox{Mor}_{e}(\C,Y)$.

Now we claim: if $H$ is a general hyperplane,
every degree $e$ effective divisors $E \leq D_H$
is a fibre of some morphism $\alpha \in \mathcal{M}$.

To prove the claim, we use the remark above that there is an open dense
subset $A \subset \PP^{r*}$ such that for every hyperplane
$H \in A$ we can find
$L$ in $W$  with $L \subset H$. We may also assume $A$ is contained in the open
set $U \!\subset\! \PP^{r*}$ of hyperplanes transverse to $C$.
By Corollary~\ref{irred} the set
$$
I_e(A)= \{(E,H) \in \C_e \times U : \; E \leq D_H, \; H \in A \}.
$$
is irreducible.
Let $Z \!\subset\! I_e (A)$ be the subset consisting of those $(E, H)$
for which the divisor $E$ is a fibre of some morphism $\alpha \in \mathcal{M}$.
The projection $pr_2: Z \ra A$ is surjective: indeed,
for every $H \in A$ we can pick $L$ in $W$ so that $L \!\subset\! H$, and
$\pi_L$ factors through some $\alpha \in \mathcal{M}$, hence
$D_H = \pi_L^{-1} (H)$ is union of  fibres of $\alpha$.
In particular,  $Z$ and  $I_e (A)$ both have dimension $r$, hence $Z$ contains
an open dense subset $Z_0$ of $I_e (A)$. Since the projection $pr_2:I_e (A) \ra A$
is surjective and generically finite, we can find an open dense subset $A_0 \!\subset\! A$
such that $pr_2^{-1} (A_0) \!\subset\! Z_0 \!\subset\! Z$, that is, for  $H \in A_0$,
every degree $e$ effective divisors $E \leq D_H$
is a fibre of some morphism $\alpha \in \mathcal{M}$, and the claim is proven.

Next we show $Y$ is rational. By way of contradiction, suppose $g(Y) \geq 2$.
Then $\mbox{Mor}_{e}(\C,Y)$ is zero dimensional, so $\mathcal{M}$ consists of a unique
morphism $\alpha$. But then not every degree $e$ subdivisor of $D_H$ can be a fibre of
$\alpha$, hence $g(Y) \leq 1$. On the other hand, if $g(Y)=1$, then two morphisms
in $\mathcal{M}$ differ only by a translation of the elliptic curve $Y$.
In particular, the morphisms in $\mathcal{M}$ have the same fibres, and as above
this is impossible. Thus $Y$ must be rational.

Finally we show $C$ is rational. We repeat the monodromy argument above,
replacing $I_e (A)$ with the set
$$
I_{e,e} (A)= \{(E_1,E_2,H) \in \C_e \times \C_e \times U : \; E_1+E_2 \leq D_H, \; H \in A \}.
$$
The conclusion is that for a general hyperplane $H$ the following
holds: for every $(E_1,E_2) \in \C_e \times \C_e$ such that $E_1+E_2 \leq D_H$
there is a subspace $L$ in $\Sigma_l$ such that $L \!\subset\! H$, the morphism $\pi_L$
factors through a morphism $\alpha: \C \ra Y \cong \PP^1$,
and $E_1$ and $E_2$ are fibres of $\alpha$.
In particular,  any two such divisors $E_1$ and $E_2$ are linearly equivalent.
Now, if $P$ and $Q$
are distinct points of $C \cap H$, we can find two divisors $F_1$ and $F_2$ of degree $e-1$
such that $P+Q+F_1+F_2 \leq D_H$. Taking first $(E_1,E_2)=(P+F_1,Q+F_2)$ and then
$(E_1,E_2)= (P+F_2,Q+F_1)$ we see that $2P$ and $2Q$ are linearly equivalent on $\C$.
Now this holds for a general hyperplane $H$ and every pair of distinct points $P$ and $Q$
in $C \cap H$, so we can  find a point $P_0 \in \C$ such that $Q-P_0$ is a
2-torsion point in the Jacobian for infinitely many points $Q$. Therefore $\C$ must be
rational.
\end{proof}

\br \label{incidence}
In the case of $\PP^3$, the second hypothesis of Theorem~\ref{trick} means
$\Sigma$ is not a Schubert cycle $\sigma(x)$ consisting of the lines
through a point $x \in \PP^3$. Indeed, given a plane $H \subset \PP^3$,
if $\sigma(H)$ denote the locus of lines contained
in $H$, it is well known that every irreducible
surface $\Sigma \!\subset\! \mathbb{G} (1, \PP^{3})$ intersects $\sigma(H)$
unless $\Sigma = \sigma(x)$ for some $x \in \PP^3$.
\er

\bc \label{decomp3}
Let $C\!\subset\!\PP^3$ be a reduced irreducible nondegenerate curve.
Suppose $C$ is not rational, and $\Sigma \!\subset\! \mathbb{G} (1, \PP^3)$
is an irreducible closed surface
whose general line does not meet $C$ and is  decomposable. Then
$$\Sigma = \sigma(x)=\{L \in \mathbb{G} (1, \PP^3): \:\: x \in L \},$$
where $x \in \PP^3 \!\smallsetminus\! C$ is one of the finitely  many non-birational
points of $C$.
\ec
\begin{proof}
By Theorem~\ref{trick} and the previous remark, we know any such surface is
a Schubert cycle $\sigma(x)$ where $x$ is a point of $\PP^3 \!\smallsetminus\! C$. If projection
from $x$ were birational onto its image, then as in the proof of Theorem~\ref{codimension}
all but finitely many lines in $\sigma(x)$ would be uniform, hence indecomposable.
\end{proof}

\section{Surfaces of non-uniform lines for curves in $\PP^3$} \label{five}
In this paragraph we consider the case of {\it smooth} irreducible curves
$C\!\subset\!\PP^3$. As in Section~\ref{three} we let
$\mathcal{N}(C)$ denote the locus of lines that do not meet $C$ and are not uniform.
and we begin the classification of pairs
$(C, \Sigma)$ where $\Sigma$ is
an irreducible surface in $\mathbb{G} (1, \PP^{3})$ such that the
general line in $\Sigma$ belongs to  $\mathcal{N}(C)$.

By Theorem~\ref{smooth-two} $\mathcal{N}(C)$
has dimension at most two, and we now wish to classify pairs
$(C, \Sigma)$ where $\Sigma$ is
an irreducible surface in $\mathbb{G} (1, \PP^{3})$ such that the
general line in $\Sigma$ belongs to  $\mathcal{N}(C)$.

First of all, as we observed in Remark~\ref{rmk},
if $x \in \PP^3 \smallsetminus C$ is one of the (finitely
many) points non birational points for $C$, then every line through $x$ is decomposable,
hence non-uniform. Thus, if $\sigma(x) \cong \PP^2$ denotes the Schubert cycle
of lines through $x$,  for each
such point $x$ the general line in $\sigma(x)$
belongs to $\mathcal{N}(C)$.
In fact we will show that, except $C$ is a rational curve of degree three, four or six,
these are the only irreducible surfaces in $\mathcal{N}(C)$. To be more precise,
we make the following definition.

\bd
Let $C$ be a nondegenerate irreducible curve in $\pso$.
We say $C$ {\it has special monodromy} if there exists an irreducible Zariski closed
surface $\Sigma \!\subset\! \mathbb{G}(1, \PP^3)$ that is not a Schubert cycle $\sigma (x)$
of lines through a point $x \in \pso$, with the property:
the general line $L \in \Sigma$ does not meet $C$, and
the projection $\pi_L$ factors as
$$C \stackrel{\alpha}{\ra} \PP^1 \stackrel{\beta}{\ra} \PP^1$$
where $\beta$ has degree $2$.
\ed
We will show in Section~\ref{six} that a curve with special monodromy is rational
of degree either $4$ or $6$. In this section we prove:

\bt \label{space}
Let $C$ be a smooth nondegenerate irreducible curve in $\pso$.
Suppose $\Sigma \!\subset\! \mathbb{G}(1, \PP^3)$ is an irreducible Zariski closed surface
such that the general line $L$ in $\Sigma$ does not meet $C$ and is not uniform.
Then one of the following three possibilities holds:
\ben
\item either there exists a non birational point $x \in \pso$, and
$\Sigma=\sigma(x)$ is the cycle of lines through $x$; or
\item
$C$ is a twisted cubic curve, and the general line in
$\Sigma$ is the intersection of two osculating planes to $C$; or
\item
$C$ is rational and has special monodromy.
\een
\et
\br
The theorem is true more generally for an irreducible nondegenerate curve
$C\!\subset \! \pso$ for which the branch divisor $B$ of
Proposition~\ref{branch-divisor} is reduced.
In particular, it
holds for curves with only nodes or simple cusps. The proof is the same as below except
$C^*$ has to be replaced by $B$.
\er

In the proof of the Theorem~\ref{space} we will need the following two lemmas.
The first one is  elementary and well known:
\bl \label{2branch}
Let $X$ be a smooth irreducible projective curve.
Suppose $\pi: X \ra \PP^1$ is a nonconstant morphism of degree $d \geq 2$. Then
$\pi$ ramifies over at least two distinct points $y_1$ and $y_2$ of $\PP^1$.
If there are no other
branch points, then $X$ is rational and up to a choice of coordinates $\pi$
is  the map $z \mapsto z^d$.
\el
The second lemma we need can be proven using a result due to Strano:
\bl[Strano] \label{Strano}
Let $X \!\subset\! \PP^3$ be an arbitrary reduced curve.
If $\Sigma \!\subset\! \mathbb{G} (1, \PP^{3})$ is an irreducible Zariski closed surface
parametrizing lines that are trisecant to $X$ (meaning $\deg (X.L) \geq 3$),
then one of the following holds:
\begin{enumerate}
\item  either $\Sigma = \sigma(x)$, where $x$ is a point in $\PP^3$ such that
$\deg_{x} (X.L) \geq 3$ for every line $L$ through $x$; or
\item $\Sigma = \sigma(H)$,
where $H$ is a plane in $ \PP^3$ containing a subcurve $P$ of $X$
of  degree at least $3$, and $\sigma(H)$ is the locus of lines
contained in $H$.
\end{enumerate}
\el
\begin{proof}
If the projection
$$\rho: J_\Sigma = \{(H,L) \in \PP^{3*} \times \Sigma: \; H \supset L \} \rightarrow \PP^{3*}$$
is not dominant, then by Remark~\ref{incidence} there is a point $x \in \PP^3$ such that
$\Sigma = \sigma(x)$, and then every line through $x$ meets $X$ in a scheme of length at least
three. This is possible only if $\deg_{x} (X.L) \geq 3$ for every line $L \in \sigma(x)$ .

Suppose now $\rho$ is dominant. This means that the general plane
contains a trisecant line $L \in \Sigma$. By \cite{Strano1}, Theorem 9,
$X$ contains a unique subcurve $P=P(L)$ whose general plane
section is $L.X$. Since $\deg (L.X ) \geq 3$, $P$ is a plane curve
of degree at least $3$. Let $H(P)$ be the linear span of $P$. Then
$\sigma(H(P))$ is a surface of trisecant lines for $X$. Since $X$
contains finitely many subcurves and the general line $L$ of $\Sigma$ is
contained in $\sigma(H(P(L)))$, $\Sigma$ must be one of the Schubert
cycles $\sigma(H(P(L)))$.
\end{proof}

\begin{proof}[Proof of Theorem~\ref{space}]
Suppose $L$ is a non-uniform line that does not meet $C$.
Then - see Remark~\ref{smooth-two} - the line $\PP^1_L \subseteq \PP^{3*}$
dual to $L$ meets the singular locus $C^{*}_{sing}$ of the dual surface $C^*$ in at
least two distinct points.

Suppose first the dual of a
general line in $\Sigma$ meets $C^{*}_{sing}$ in more than two distinct points.
Let $X$ denote the curve $C^{*}_{sing}$ taken with its reduced scheme structure,
and apply Lemma~\ref{Strano} to $X$
and to the surface $\Sigma' = \{\PP^1_L: L \in \Sigma \}$. The conclusion is that either
there is a point $H^*$ such that $\Sigma' =\sigma (H^*)$, or there
is a point $x \in \PP^3$ such that $\Sigma = \sigma (x)$.
We cannot have $\Sigma' =\sigma (H^*)$
because the general line through $H^*$ in $\PP^{3*}$ meets the curve $X$ only at
$H^*$, while we are assuming that a general line in $\Sigma'$ meets $X$ in more than two points.
Thus $\Sigma=\sigma(x)$. Furthermore, as in the proof of~\ref{decomp3},
the point $x$ is not birational.

Suppose now a general line $\PP^1_L$ in $\Sigma'$ meets $X$ - the singular locus of $C^*$-
in exactly two points $H_1^*$ and $H_2^*$.
If $\PP^1_L$  meets the dual surface $C^*$ only at $H_1^*$ and $H_2^*$,
the branch divisor of $\pi_L$ is supported at $H_1^*$ and $H_2^*$. Then,
by Lemma~\ref{2branch}, the curve
$C$ is rational
and $\pi_L: C \ra \PP^1$ is given in some coordinates by $z \mapsto z^d$.
Thus  $L$ is the intersection
of two planes $H_1$ and $H_2$ each of them
meeting $C$ at a single point with multiplicity $d$. As $L$ varies in a surface, there must
be infinitely many such osculating planes, and this is possible (since we are in char. zero)
only if $d=3$. Thus $C$ is a twisted cubic curve, and
$\Sigma$ consists of those lines that are intersection of two osculating planes to $C$.

We are left with the case when,  for a general $L$ in $\Sigma$, the line $\PP^1_L$
meets $C^{*}_{sing}$ in exactly two points $H_1^*$ and $H_2^*$, but it also
contains a smooth point $H^*_3 \in C^*$.
Then, by~\ref{branch-divisor} and~\ref{scambi},  the monodromy group of
$\pi_L: \C \ra \PP^1_L$ contains a
transposition. Therefore by~\ref{decompose}
the map $\pi_L$ factors nontrivially
as
$\C \stackrel{\alpha}{\ra} Y \stackrel{\beta}{\ra} \PP^1$.
If the map $\beta$ ramified over more than two distinct points
of $\PP^1_L$, these points would all be singular points of $C^*$
by~\ref{branch-divisor} because $\deg (\alpha) \geq 2$, contradicting the assumption
that $\PP^1_L$ meet $C^{*}_{sing}$ in exactly two distinct points.
Now Lemma~\ref{2branch} implies $Y \cong \PP^1$ is rational, and up to a choice of
coordinates
$\beta: \PP^1 \ra \PP^1$ is the map  $z \mapsto z^e$ where $e = \deg (\beta)  \geq 2$.
We cannot have $e \geq 4$ because as above this would imply the existence
of infinitely many hyperosculating planes to $C$. We claim the case $e=3$ is also impossible.
Since $\deg (\alpha) \geq 2$, if we had $e=3$, we would have either infinitely many
planes that osculate $C$ at more than one point, or infinitely many hyperosculating
planes. But $C$ has only finitely many hyperosculating planes, while
the first alternative contradicts the fact that the
map sending a point $P \in \C$ to the osculating plane at $P$
is birational.
Thus $e=\deg(\beta)=2$.
Finally, if $\Sigma$ is not a cycle $\sigma(x)$, the curve $C$ must be rational
by Corollary~\ref{decomp3}.
\end{proof}

\section{Rational space curves admitting surfaces of decomposable lines}\label{six}
In this section we classify smooth curves in $\pso$ having special
monodromy. The main result is that such a curve is  either a
rational quartic or a rational sextic.

We begin with some preliminary remarks. Suppose $C\!\subset\!\PP^3$ is a degree $2d \geq 4$
smooth rational curve with special monodromy. Recall this means
there exists an irreducible two dimensional family
of lines $\Sigma $, not contained in a Schubert cycle $\sigma (x)$, with the property:
every line $L \in \Sigma$ in the family does not meet $C$, and
the projection $\pi_L$ factors as
$$C \stackrel{\alpha_L}{\ra} \PP^1_L \stackrel{\beta_L}{\ra} \PP^1_L$$
where $\beta_L$ has degree $2$.

As above, we denote by $V \subset \mbox{H}^0(\pone, \coo_{\pone}
(2d))$ the $4$-dimensional subspace defining the embedding
$C\!\subset\!\PP^3=\PP (V^*)$. Since the morphisms $\beta_L$ have
degree two, the vector space $V$ contains a lot of forms $g^2$
with $g \in \mbox{H}^0(\pone, \coo_{\pone} (d))$. To make this
precise, from now on we let $\PP_m = \PP( \mbox{H}^0(\pone,
\coo_{\pone} (m)))$ and look at the embedding $q: \PP_{d} \ra
\PP_{2d}$ that sends (the class of) a degree $d$ form to (the
class of) its square. We will refer to $q$ or to its image
$X \!\subset\! \PP_{2d}$ as the {\it quadratic Veronese}.

\bp \label{one} Let $C\!\subset\!\PP(V^*)$ be a smooth rational
curve of degree $2d \geq 4$ with special monodromy. Then the
intersection $X \cap \PP(V) \subset \PP_{2d}$ has dimension at
least one. \ep
\begin{proof}
By definition of special monodromy, the curve $C$ admits
a two dimensional family $\Sigma \subset \mathbb{G} (1, \PP (V^*))$ of decomposable
lines such that,  for every $L \in \Sigma$,
the projection $\pi_L$ factors through a degree two morphism
$\beta_L: \PP^1 \ra \PP^1_L$. Up to a choice of coordinates,
$\beta_L$ is the map $z \mapsto z^2$, hence the line $\PP^1_L
\subset \PP(V) \subset \PP_{2d}$ intersects $X$ in at least two
points. Thus the surface
$$\Sigma' = \{ \PP^1_L \in  \mathbb{G} (1,\PP(V)): L \in \Sigma \} $$
is contained in the set of lines secant to $X  \cap \PP(V)$.
Hence $X \cap \PP(V)$ is at least one dimensional.
\end{proof}

Thus, if $C$ has special monodromy, we can find a curve $Y \subset
\PP_{d}$ such that $q(Y)$ lies in the three dimensional subspace
$\PP(V)$ of $\PP_{2d}$, and one expects there are very few such
curves. In fact, we show below that $Y$ must either be a conic in
$\PP_2$ or a twisted cubic in $\PP_3$.

Our first remark is that, since $C$ is smooth, we can exclude the
case $Y$ is a line: \bl \label{cone} Let $C\!\subset\!\PP(V^*)$ be
a smooth rational curve of degree $2d \geq 4$. For any line $Y
\!\subset\! \PP_d$ the conic $q(Y)$ is not contained in $\PP(V)$.
\el
\begin{proof}
The linear span of the conic $q(Y)$ is a plane $\PP(W) \!\subset\! \PP(V)$,
and the linear series $\PP(W)$ defines a projection
$C \ra \PP^2=\PP(W^*)$ whose image is a conic.
Thus $C$ is contained in a quadric cone.
But on the quadric cone there are no smooth rational curves of degree greater than $3$,
for example by~\cite[Exercise 2.9 p. 384]{AG}.
\end{proof}

We can also say that $q(Y)$ spans $\PP(V)$. One way to prove this
is to observe: \bl \label{quattro} The quadratic Veronese
$X=q(\PP_d) \!\subset\! \PP_{2d}$ contains no four collinear
points. \el
\begin{proof}
Fix a line  $M \!\subset\! \PP_{2d} $, and think of it as a linear
series $g^1_{2d}$ on $\PP^1$. Suppose $M$ contains two distinct
points of $X$. Then $M$ corresponds to a subspace $\langle g^2,h^2
\rangle$ of $\mbox{H}^0(\pone, \coo_{\pone} (d))$. It follows
that, if $M$ has a base point $p \in \PP^1$, then in fact $2p$ is
in the base locus of $M$. We can then remove $2p$ from the base
locus and deduce the statement from the case $d\!-\!1$ by
induction. Therefore we can assume $M$ is base point free and
defines a morphism $\pi: \PP^1 \ra \PP^1$ of degree $2d$. By
Riemann-Hurwitz the ramification divisor has degree
$$
\deg (R_{\pi}) = 4d-2.
$$
Now, if the linear series $M$ contains a divisor $2D$, then $D
\!\subset\! R_{\pi}$. Since $D$ has degree $d$, there are at most
$3$ such divisors in $M$, that is, the support of $X \cap M$
consists of at most $3$ points (by contrast, $M$ can be tangent to
$X$ at one or more points of intersection, and the scheme
theoretic intersection $X \cap M$ may well have length~$4$).
\end{proof}

\bc \label{span} Let $C\!\subset\!\PP(V^*)$ be a smooth rational
curve of degree $2d \geq 4$. Suppose $Y \subset \PP_d$ is a curve
such that $q(Y) \subset \PP(V)$. Then $\PP(V)$ is the linear span
of $q(Y)$. \ec
\begin{proof}
Let $T=q(Y)$, and note $T$ is a curve in $X \cap \PP(V)$. Since
$X=q(\PP_d)$ is an embedding of $\PP_d$ defined by a linear system
of quadrics, every curve in $X$ - and in particular $T$ - has even
degree. Now $T$ is not a conic by~Lemma~\ref{cone}, and it cannot
be a plane curve of degree $\geq 4$ because $X$ contains no $4$
collinear points. Thus $T$ spans $\PP(V)$.
\end{proof}

Before proving $Y$ must be a conic or a twisted cubic, we recall
results due as far as we know to Hopf and
Eisenbud~\cite[5.2]{eisenbud}. We denote the linear span of a
subset $S\!\subset\!\PP^N$ by the symbol $\langle S \rangle$.

\bp \label{linear} Let $q: \PP_d \ra \PP_{2d}$ be the quadratic
Veronese. \ben
\item (Hopf)
If $H\!\subset\! \PP_d$ is a linear subspace of dimension $r$,
then
$$ \dim \, \langle q(H) \rangle \geq 2r.$$
\item(Eisenbud)
Equality holds if and only if there exist a two dimensional
subspace $W\!\subset\! \PP_a$ and a form $h \in \PP_{d-ra}$ such
that $H = \PP(h \mbox{S}^r (W))$ where $\mbox{S}^r (W)$ is the
image of $\mbox{\em Sym}^r (W)$ in the space of forms of degree
$ra$ - that is, there exist two linearly independent forms $F$ and
$G$ of degree $a$ such that $H = \PP (\langle hF^r, hF^{r-1} G,
\ldots, hG^r \rangle)$.
\een \ep

\bp \label{final} Suppose $Y \!\subset\! \PP_d$ is a reduced
irreducible curve of degree $e$ such that $q(Y)$ spans a three
dimensional linear subspace $\PP(V)$ of $\PP_{2d}$. Suppose
further the linear series $\PP(V)$ is base point free and defines
an embedding $f: \PP^1 \ra \PP(V^*)$ with image $C$. Then either
$d=e=2$, or $d=e=3$ and $Y$ is a twisted cubic. In particular, $C$
is either a rational quartic or a rational sextic. \ep
\begin{proof}
Given a point $s \in \PP^1$, we denote by $H(s) \!\subset\! \PP_d$
the hyperplane of forms vanishing on $s$. Note that by assumption
$\PP(V)$ is base point free, hence $Y$  is not contained in a
hyperplane $H(s)$. We first prove the following claim:

\vspace{0.2cm}
\begin{em}
For general $s \in \PP^1$, the intersection
of $H(s)$ and $Y$ is transverse, i.e., consists of  $e$ distinct points.
\end{em}
\vspace{0.2cm}

Since $Y$ has only finitely many singular points and each of them is contained
in only finitely many hyperplanes of the form $H(s)$, the general hyperplane
$H(s)$ does not meet the singular locus of $Y$. Thus if
the claim were false, for every $s$ in $\PP^1$ the hyperplane $H(s)$ would be tangent to
$Y$. We contend this implies $Y \!\subset\! \Delta$ where $\Delta$ is the
hypersurface of polynomials having a double root.

Analytic proof: consider the normalization morphism $g: \widetilde{Y} \ra Y$
and the correspondence
$$
J= \{ (x,s) \in \widetilde{Y} \times \PP^1: \,\, T_x Y \!\subset\! H(s) \; \}.
$$
We are assuming the second projection $J \ra \PP^1$ is surjective. Since
$J$ is one dimensional, we can find near a general point $s \in \PP^1$
a local analytic section $s \mapsto (x(s),s)$ of $p_2: J \ra \PP^1$. Choose
local coordinates $x$ on $\widetilde{Y}$ and $s$ on $\PP^1$ and write $g(x,s) \in \mathbb{C}$
for the value of the polynomial $g(x)$ at $s$. The condition $T_{x(s)} Y \!\subset\! H(s)$
means $g(x(s),s) = \partial_x g(x(s),s) = 0$. Therefore
$$
0 = \frac{d}{ds} \, g (x(s),s) = \partial_x g(x(s),s) x'(s) + \partial_s g(x(s),s) =
\partial_s g(x(s),s)
$$
so that $s$ is a double root of the polynomial $g(x(s))$. It follows that
$Y \!\subset\! \Delta$.

Geometric proof: to say  for general $s$ the hyperplane $H(s)$ is tangent
to $Y$ is to say the hypersurface $Y^* \!\subset\! \PP^{d*}$  dual to $Y$
contains the rational normal curve $$\Gamma = \{ H(s) : \, s \in \PP^1\}.$$

Since $\Gamma \!\subset\! Y^*$, a tangent line $L$ to $\Gamma$ is contained in
some hyperplane $\Pi_L $ tangent to $Y^*$, i.e., belonging to $Y^{**}$.
Then $\Pi_L \in Y^{**} \cap \Gamma^* = Y \cap \Delta$, where the last
equality follows from the biduality theorem and the fact that $\Gamma$ is the dual
of $\Delta$. Now $\Pi_L$ has to vary in a one
dimensional family because $\Gamma$ is nondegenerate. Hence $Y \cap \Delta$
is one dimensional, and therefore must coincide with $Y$.

We have shown $Y \!\subset\! \Delta$, that is, every form in $Y$ has at least one
root of multiplicity at least $2$. It follows the linear series $\PP(V)$ spanned
by $q(Y)$ contains infinitely many divisors with multiplicity at least $4$ at one
point. But, since $\PP(V)$ is three dimensional, this implies the linear series has
infinitely many inflectionary points, a contradiction. Therefore the claim holds:
for $s$ general in $\PP^1$ the hyperplane $H(s)$ meets $Y$ in $e = \deg (Y)$ distinct
points.

We now want to apply Riemann-Hurwitz to a general tangential projection of $C$
and use the claim to prove that, if $e$ is not small, then the ramification divisor
is too large. By a theorem of Kaji~\cite{kaji}, for $s$ general  in $\PP^1$,
the tangent line $T_{f(s)} C$
meets $C$ only at $f(s)$ and with multiplicity two.
Hence projection from $T_{f(s)} C$ is a degree $2d-2$ morphism
$\pi: C = \PP^1 \ra \PP^1$
In other words, $\pi$ is the morphism defined by the base point free $g^1_{2d-2}$:
$$
\{D - 2s : \: D \in |V|, D \geq 2s    \}.
$$
Now let $g$ be a form in $Y \cap H(s)$: then $2s$ is a double root
of  $g^2 \in \PP(V)$, hence our $g^1_{2d-2}$ contains the divisor
$$
(g^2)_0 - 2s = 2s_1 + \cdots + 2 s_{d-1}
$$
and we see  $s_1 + \cdots + s_{d-1} \leq R_{\pi}$ where $R_{\pi}$ is the ramification
divisor of $\pi$. Since $Y \cap H(s)$ consists of $e$ distinct points, we conclude
$\deg R_{\pi} \geq e (d-1)$. By Riemann-Hurwitz then
$$
-2 = -2 (2d-2) + \deg (R_{\pi}) \geq -2 (2d-2) + e (d-1)
$$
so that $$ e \leq 4 - \frac{2}{d-1}.$$
Since $e>1$ (otherwise the span of $q(Y)$ would be a $\PP^2$), we must have
either $e=3$ and $d \geq 3$ or $e=2$ and $d \geq 2$.

Suppose first $e=2$. Then $Y$ is a smooth conic spanning a
2-dimensional linear subspace $H \subseteq \PP_d$. Since $q_{| H}:
H \ra \PP_{2d}$ is defined by a linear system of quadrics, the
linear space $\langle q(H) \rangle$ has dimension at most $5$, and
by Hopf at least $4$. If it has dimension $5$, then $q_{| H}$ is
the Veronese embedding of $H = \PP^2$, and then $q$ maps any conic
of $\PP^2$ to a hyperplane section of the Veronese surface, and so
$q(Y)$ spans a $\PP^4$, contrary to our assumptions. Hence
$\langle q(H) \rangle$ is a $\PP^4$.

Now by Eisenbud's result~\ref{linear}(2) we have
$$\langle q(H) \rangle =
\PP( \langle h^2 F^4, h^2 F^3G, h^2 F^2G^2, h^2 FG^3, h^2 G^4 \rangle).$$
Thus the linear series $\langle q(H) \rangle $ defines a morphism $\PP^1 \ra \PP^4$ whose
image is a normal rational quartic $T$, and the linear subseries $\PP(V)$ maps $\PP^1$ to a
projection $C$ of $T$ in $\PP(V^*)$. Therefore $\deg (C) \leq 4$, that is, $d = 2$.

Finally suppose $e=3$ and $d \geq 3$. First observe $Y$ is a
twisted cubic curve. Otherwise $Y$ would be a plane cubic and
would have infinitely many trisecant lines. This is impossible
since $Y \!\subset\! q^{-1} (\PP(V))$ and $q^{-1} (\PP(V))$ is cut
out by quadrics. Thus $Y$ is a twisted cubic curve spanning a
three dimensional subspace $H \!\subset\! \PP_d$. Let $ \PP^s = \;
\langle q(H) \rangle $. By hypothesis $q(Y)$ spans a $\PP^3$ in
$\PP^s$, so $Y$ is contained in $s-3$ linearly independent
quadrics. But  $Y$ is a twisted cubic curve, thus $h^0 (H,
\ideal{Y,H} (2)) = 3$ and we conclude $s \leq 6$.  By Hopf, we
must have $s=6$. As in the case $e=2$, we can now use Eisenbud's
result~\ref{linear}(2) to conclude $d=3$ and $C$ is a rational
normal sextic curve.
\end{proof}

\br
The hypothesis $f: \PP^1 \ra \PP(V^*)$ is an embedding - rather than
birational onto its image - is used only to
apply the theorem of Kaji (a general tangent to $C$ intersects $C$ only at the
point of tangency). But one can avoid using this fact analyzing what happens
if the projection from a general tangent line has degree less than $2d-2$.
\er

\br \label{whichsurfaces}
In fact, in both the cases $d=2$ and $d=3$, one sees $q(Y) = X \cap \PP(V)$
for example using the fact $q^{-1} (X \cap \PP(V))$ is cut out by the right
number of quadrics. Thus $Y$ is unique, and the surface $\Sigma$ of decomposable
lines for $C$ is determined as the dual of the families of lines in $\PP(V)$
that are secant to the rational curve $q(Y)$.
\er

\bex
The proposition is false in positive characteristic.
An example in char. 7 of a normal rational quartic $Y \!\subset\! \PP_4$
such that $q(Y)$ spans only a $\PP^3$ in $\PP_8$ is given by the parametric equation
$$
g(t) = 1 - 2 t X - 2 t^2 X^2 - 4 t^3 X^3 + 4 t^4 X^4
$$
Indeed $q(Y)$ is parametrized by
$$
g(t)^2 =
1 - 4 t X + 28 t^4 X^4 - 32 t^7 X^7 + 16 t^8 X^8.
$$
\eex

\bt \label{end}
Suppose $C$ is a smooth rational curve having special monodromy. Then
$\deg C$ is either $4$ or $6$. Furthermore, every rational quartic
has special monodromy; the general rational sextic does not have special monodromy,
but there are rational sextics with special monodromy.
\et
\begin{proof}
Since $C$ has special monodromy, its degree is an even number $2d
\geq 4$ and by Proposition~\ref{one} there is an irreducible curve
$Y \!\subset\! \PP_d$ such that $q(Y) \!\subset\! \PP(V)$. By
Corollary~\ref{span} the curve $q(Y)$ spans $\PP(V)$. Now
Proposition~\ref{final} implies $\deg (C)$ is either $4$ or $6$.

To see that every rational quartic has special monodromy, we observe the converse
of Proposition~\ref{one} holds: if $C\!\subset\!\PP(V^*)$ is a smooth rational curve of
degree $2d \geq 4$ and $\dim X \cap \PP(V) \, > 0$, then $C$ has special monodromy.
Indeed, suppose there is an irreducible curve $T \subset X \cap \PP(V)$.
By Corollary~\ref{span}, the curve $T$ spans $\PP(V)$, that is, it is not contained
in a plane. Therefore the secant lines to $T$ are not all
contained in one plane. Now let $\Sigma' \!\subset\!  \mathbb{G}
(1,\PP(V))$  be the surface of secant lines to  $T$, and let
$\Sigma$ be the corresponding subvariety of $\mathbb{G} (1,\PP(V^*))$.
Since the secant lines to $T$ are not
all contained in one plane, $\Sigma$ is not contained in a
Schubert cycle $\sigma(x)$. Since $\PP(V)$ is spanned by $T$, a general
secant line $\PP^1_L$ to $T$ is base point free, so that its dual line
$L \in \Sigma$ does not meet $C$. Furthermore, by construction a general
line in $\Sigma$ is decomposable as prescribed by special
monodromy.  Therefore $C$ has special monodromy.

Now suppose $C$ is a rational quartic, that is, $d=2$.
Then  $X$ is a Veronese surface in $\PP_4$, and every
hyperplane $\PP(V)$ in $\PP_4$ cuts a curve on $X$. Thus $C$ has special monodromy.

When $d =3$, the quadratic Veronese $X$ has codimension $3$ in $\PP_6$, therefore a
general $\PP^3$ in $\PP_{6}$ does not cut a curve on $X$, and the general
rational sextic does not have special monodromy. For an example of a sextic
with special monodromy, consider - cf. \cite{bm} -
the twisted cubic curve $Y \!\subset\! \PP^3 \cong \mathbb{C}[X]_{\leq 3}$,
image of the embedding $g: \PP^1 \ra \PP^3$ defined by
$$
g(t)= X + tX^3 + t^2 - t^3X^2.
$$
One immediately checks $q(Y)$ spans a three dimensional subspace $\PP(V)$ in $\PP_6$,
and the linear series $|V|$ defines an embedding of $\PP^1$ in $\PP(V^*)$: the image
$C$ of this embedding is a rational sextic with special monodromy. This example
can be constructed noting $q(Y)$, as every sextic in $\PP^3$, must have quadrisecant lines,
and these must be degenerate because of Lemma~\ref{quattro}. Thus one looks
for $Y$ such that the line joining $g(0)^2=X^2$ and $g(\infty)^2 = X^4$ is tangent to $q(Y)$
both at $g(0)^2$ and at $g(\infty)^2$.
\end{proof}

\section{Examples of curves admitting a one dimensional family of non~uniform~lines}
\label{seven}
We give examples of space curves having a one dimensional
  family of non uniform lines.
   The monodromy of all these lines
  is contained in the alternating group.

   We recall
  the construction given in~\cite{pi}.  Let $X$ be any smooth curve of genus $g$ and
   $L$ be a spin line bundle on $X,$ i.e.
   $L^{2}=\omega_{X}.$ We fix a point $P$ in $X$ - in fact one could
    take any effective divisor $D>0$.
   Let $V_{m}= H^{0}(X,L(mP))$ and
   $W_{m}=H^{0}(X,\omega_{X}(2mP))$.
   Given a two plane $\Pi \subset V_{m}$,  we let
   $S^{2}(\Pi)\subset W_{m}$ be the $3-$dimensional space image of the natural map
   $\Pi \otimes \Pi \to W_{m}.$ If $(s,t)$ is a base of $\Pi$,
   we define a basis of $S^{2}(\Pi)$ setting $ \omega_{1}=s^2-t^2,$
   $\omega_{2}=i(s^{2}+t^{2})$ and
   $\omega_{3}= 2st.$ In~\cite{pi}, a plane $\Pi$ is called {\em minimal} if all the elements
   of $S^{2}(\Pi)$ are exact meromorphic forms.
  Then there exist $F_{i}\in H^{0}(X,\mathcal O_{X}(2m-1)P))$ such that
    $dF_{i}= \omega_{i}.$ The name minimal and the above choice of the basis
     come from the fact the real parts of the $F_{i}$ define a map $G: X-P \ra \mathbb{R}^3$
     $$G(x)= Re(F_{1}(x),F_{2}(x),F_{3}(x)),$$
     which gives  a minimal surface in the Euclidean space.
     Associated to a minimal plane $\Pi$
    we  consider the space:
   $$ W =\{ k\in H^{0}(X,\mathcal O_{X}(2m-1)P) : dk \in
S^{2}(\Pi)\}.$$ Note that $W$ is generated
   by the constants and the $F_{i}$. If moreover the linear system associated
    to $\Pi$ is base point free, then
   $W$ is also base point free and it defines an immersion
   $f: X\to \mathbb P^{3}$ setting $f(x)=(1,F_1,F_2,F_{3}).$
  For suitable choices (for instance if 2m-1 is a prime number)
   $f(X) $ is a curve of degree 2m-1 and hence $f$ is birational onto
    its image.

    The existence, for any $X$ and $P$,
     of base point free minimal planes $\Pi\subset V_{m}$ is proven
     in~\cite[Proposition 5.8 p. 355]{pi}
     for  $m>31g+21$.
   On the other hand, if $\sigma \in \Pi, \sigma\neq 0,$ let
   $ L_{\sigma}\subset W$ be the two
   dimensional subspace defined by the equation
    $dF =\lambda\sigma^{2}$, $\lambda \in \mathbb C$. Note
   $L_{\sigma} = \langle 1,F_{\sigma} \rangle $ where $dF_{\sigma}=\sigma^{2}.$
   Now $ L_{\sigma} $ defines a non uniform line
   for $f(X)$. In fact  the monodromy of
   $F_{\sigma}$  is contained in the alternating group. This follows because all
    the ramifications
     $F_{\sigma}$ are odd. The map $\sigma \to L_{\sigma}$
   defines a $\mathbb P^1=\mathbb P(\Pi)$ family of non uniform lines
for $C=f(X)$.
   We remark that all the lines belong to a fixed  plane.

\end{document}